\documentstyle[amssymb,12pt]{article}
\def\ker{{\rm ker}}
\def\cali{{\cal I}}
\def\calp{{\cal P}}
\def\cala{{\cal A}}
\def\calh{{\cal H}}
\def\oplusinf{\mathop{\oplus}}
\def\cals{{\cal S}}
\def\calc{{\cal C}}
\def\calf{{\cal F}}
\def\k{\ \Bbb K\ }
\def\Bbbone{\mbox{\rm 1\hspace {-.6em} l}}

\newcommand{\vir}{\raisebox{0.75mm}{,}}
\newtheorem{theorem}{THEOREM}
\newtheorem{lemma} {LEMMA}

\begin{document}

\baselineskip=0.7cm
\begin{center}
{\Large\bf BASIC COHOMOLOGY \\
OF ASSOCIATIVE ALGEBRAS}
\end{center}
\vspace{0.75cm}

\begin{center}
{\large Michel DUBOIS-VIOLETTE and Thierry MASSON}
\end{center}
\vspace{0.2cm}
\begin{center}
{\small\it Laboratoire de Physique Th\'eorique\footnote{Laboratoire associ\'e
au
C.N.R.S.}\\ B\^atiment 211, Universit\'e Paris XI\\
F-91405 ORSAY Cedex}\\

E-mail: flad@qcd.th.u-psud.fr and masson@qcd.th.u-psud.fr\\

\today
\end{center}

\vspace {1cm}
\begin{abstract}
We define a new cohomology for associative algebras which we compute for
algebras with units.
\end{abstract}
\vspace{1cm}

L.P.T.H.E.-ORSAY 94/10
\newpage
\section{Introduction: Definition of the basic cohomology of an
associative algebra}

Let $\cal A$ be an associative algebra over $\k=\Bbb R$ or $\Bbb C$ and
let ${\cal A}_{Lie}$ be the underlying Lie algebra (with the commutator
as Lie bracket). For each integer $n\in\Bbb N$, let $C^n({\cal A})$ be
the vector space of $n$-linear forms on $\cal A$, i.e. $C^n({\cal
A})=(\cala^{{\otimes}^{n}})^{\ast}$. For $\omega\in C^n(\cala)$ and $\tau
\in C^m(\cala)$ one defines $\omega . \tau\in C^{n+m}(\cala)$ by:
$$\omega .\tau(A_1,\dots,A_{n+m})=\omega(A_1,\dots,A_n)\tau
(A_{n+1},\dots,A_{n+m}),\ \forall A_i\in \cala.$$
Equipped with this product, $C(\cala)=\displaystyle{\oplusinf_n} C^n(\cala)$
becomes an
associative graded algebra with unit $(C^0(\cala)=\k )$. One defines a
differential $d$ on $C(\cala)$ by setting for $\omega\in
C^n(\cala)$, $A_i\in \cala$
$$d\omega(A_1,\dots,A_{n+1})=\sum^n_{k=1}(-1)^k\omega(A_1,\dots,A_{k-1},A_kA_{k+1},A_{k+2},\dots,A_{n+1}).$$
Indeed, $d$ is the extension as antiderivation of $C(\cala)$ of minus
the dual of the product of $\cala$ and $d^2=0$ is then equivalent to the
associativity of the product of $\cala$. The graded differential algebra
$C(\cala)$ is together with a bimodule $\cal M$ the basic building
blocks of the Hochschild complex giving the Hochschild cohomology with value
in $\cal M$. Here we do not want to introduce bimodules like $\cal M$.
However it is well known, see below, that the cohomology of $C(\cala)$
is trivial whenever $\cala$ has a unit. Nevertheless, there are two
classical cohomologies which can be extracted from the differential
algebra $C(\cala)$, namely the Lie algebra cohomology of $\cala_{Lie}$
and the cyclic cohomology of $\cala$. In fact let
$\cals:C(\cala)\rightarrow C(\cala)$ and $\calc :C(\cala)\rightarrow
C(\cala)$ be defined by $$(\cals\omega)(A_1,\dots,A_n)=\sum
_{\pi\in \cals_{n}}\varepsilon(\pi)\omega (A_{\pi(1)},\dots,A_{\pi(n)})$$ and
$$(\cal C\omega)(A_1,\dots,A_n) =\sum_{\gamma\in{\cal C}_{n}}\varepsilon
(\gamma) \omega (A_{\gamma(1)},\dots,A_{\gamma(n)})$$
for $\omega\in C^n({\cal A})$, $A_k\in {\cal A}$ and where $\cals_n$ is the
group of permutations of $\{1,\dots,n\}$ and ${\cal C}_n$ is the
subgroup of cyclic permutations. One has $\cals\circ d=\delta \circ \cals$
where
$\delta$ is the Chevalley-Eilenberg differential so $(Im\ \cals,\delta)$ is
a differential algebra whose cohomology is the Lie algebra cohomology
 $H({\cal A}_{Lie})$ of the Lie algebra
${\cala}_{Lie}$ [3],[6],[5]. On the other hand, see Lemma 3 in [4] part II, one
has
${\cal C}\circ d = b\circ {\cal C}$ where $b$ is the Hochschild
differential of $C({\cala},\cala^\ast)$ so $(Im\ {\cal C},b)$ is a complex
whose cohomology is the cyclic cohomology $H_\lambda(\cala)$ of $\cala$
up to a shift $-1$ in degree [4], (it is worth noticing, and this is not
accidental, that the same shift occurs in the Loday-Quillen theorem [7]).\\
We want now to point out that there is another natural non-trivial
cohomology which may be extracted from the differential algebra $C({\cal
A})$. This cohomology is connected with the existence of a canonical
operation, in the sense of H. Cartan [2], [5], of the Lie algebra
$\cala_{Lie}$ in the graded differential algebra $C(\cala)$. For $A\in
\cala=\cala_{Lie}$, define  $i_A:C^n(\cala)\rightarrow C^{n-1}(\cala)$
by
 $$i_A\omega(A_1,\dots ,A_{n-1})=\sum^{n-1}_{k=0}(-1)^k \omega
(A_1,\dots ,A_k,A,A_{k+1},\dots,A_{n-1})$$
 $\forall \omega\in
C^n(\cala)$, $\forall A_i\in \cala$, for $n\geq 1$ and $i_AC^0(\cala)=0$.
For each $A\in \cala$, $i_A$ is an antiderivation of degree $-1$ of
$C(\cala)$ and one has, with $L_A=i_Ad+di_A$, $i_Ai_B+i_Bi_A=0$,
$[L_A,i_B]=i_{[A,B]}$, $[L_A,L_B]=L_{[A,B]}$ which are the relations
which characterize an operation of $\cala_{Lie}$ in $C(\cala)$. Notice
that then, for $A\in\cala$, the derivation $L_A$ of degree 0 of
$C(\cala)$ is given by
$$L_A\omega(A_1,\dots,A_n)=\sum^n_{k=1}\omega (A_1,\dots,
[A_k,A],\dots,A_n)$$ for $\omega\in C^n(\cala)$, $A_i\in \cala$. An
element $\omega\in C(\cala)$ is called {\it horizontal} if $i_A\omega=0$
for any $A\in \cala$, it is called {\it invariant} if $L_A\omega=0$ for
any $A\in \cala$ and it is called {\it basic} if it is horizontal and
invariant, i.e. if $i_A\omega=0$ and $L_A\omega=0$ for any $A\in \cala$.
The set $C_H(\cala)$ of horizontal elements of $\cala$ is a graded
subalgebra of $C(\cala)$ which is stable by the $L_A$, $A\in \cala$. The
set $C_I(\cala)$ of invariant elements of $\cala$ and the set
$C_B(\cala)$ of basic elements of $\cala$ are two graded differential
subalgebras of $C(\cala)$ $(C_B(\cala)\subset C_I(\cala))$; their
cohomologies $H_I(\cala)$ and $H_B(\cala)$ are called the invariant
cohomology and the basic cohomology of $\cala$. As already claimed, if
$\cala$ has a unit then the cohomology $H(\cala)$ of $C(\cala)$ is
trivial and it turns out that the same is true for the invariant
cohomology ; one has the following proposition.
\newtheorem{proposition}{PROPOSITION}
\begin{proposition}
If $\cala$ has a unit, then one has $H^n(\cala)=0$, $H^n_I(\cala)=0$ for
$n\geq 1$ and $H^0(\cala)=H^0_I(\cala)=\k$.
\end{proposition}

\noindent {\bf Proof.} Let $\Bbbone$ be the unit of $\cala$ and let us
define for $n\geq 1$\\
 $h:C^n(\cala)\rightarrow C^{n-1}(\cala)$ by
$h\omega(A_1,\dots,A_{n-1})=-\omega(\Bbbone,A_1,\dots,A_{n-1})$,\\ for
$\omega\in C^n(\cala)$ and $A_i\in
\cala$. One has $(dh+hd)\omega=\omega$ and $(L_Ah-hL_A)\omega=0$ for
$\omega\in C^n(\cala)$ and $A\in \cala$. It follows that $h$ is a
contracting homotopy for $C^+(\cala)=\displaystyle{\oplusinf_{n\geq
1}}C^n(\cala)$ and for $C^+_I(\cala) = \displaystyle{\oplusinf_{n\geq
1}}C^n_I(\cala)$, which proves the
result.$\square$

The basic cohomology $H_B(\cala)$ is however non-trivial. In fact it is
already non-trivial for $\cala=\k$.

\begin{proposition}
The basic cohomology $H_B(\k)$ of $\k$ is the free graded
commutative algebra with unit generated by an element of degree
two;\linebreak[4]
$H^{2k}_B(\k)=\k$, $H^{2k+1}_B(\k)=0$ and $H_B(\k)$
identifies to the algebra $\k[X^2]$ of polynomials in one
indeterminate $X^2$ of degree two, $(X^2$ being identified to a
non-vanishing element of $H^2_B(\k))$.
\end{proposition}

\noindent {\bf Proof.} $C(\k)$ can be identified to $\k[X]$ and
coincides with $C_I(\k)$ since $L_1=0$. One has $i_1=0$ on the
elements of even degrees and $i_1\not= 0$ on the non-vanishing elements
of odd degrees.\\ Therefore $C_B(\k)=\displaystyle{\oplusinf_k}
C^{2k}(\k)=\k [X^2]=H_B(\Bbb K).$ $\square$

It is worth noticing here that one has $C^1_B(\cala)=0$ and therefore
$H^1_B(\cala)=0$ for any associative $\k$-algebra $\cala$.\\
In the next section we shall compute $H_B(\cala)$ for an arbitrary
associative $\k$-algebra $\cala$ with unit.

\section{Computation of the basic cohomology of unital algebras}

In this section, $\cala$ is an associative $\k$-algebra with a unit
denoted by $\Bbbone$. Let $\cal I^n_S(\cala_{Lie})$ denote the space of
ad*-invariant homogeneous polynomials of degree $n$ on the underlying
Lie algebra $\cala_{Lie}$ of $\cala$. We shall prove the following
theorem which generalizes the proposition 2 of \S 1.

\begin{theorem}
The basic cohomology $H_B(\cala)$ of $\cala$ identifies with the algebra
$\cal I_S(\cala_{Lie})$ of invariant polynomials on the Lie algebra
$\cala_{Lie}$ where the degree $2n$ is given to the homogeneous
polynomials of degree $n$, i.e. $H^{2n}_B(\cala)\simeq \cal
I^n_S(\cala_{Lie})$ and $H^{2n+1}_B(\cala)=0$. In particular,
$H_B(\cala)$ is commutative and graded commutative.
\end{theorem}

In order to prove this theorem, we shall need some constructions used in
equivariant cohomology [1]. Let $\calp^{m,n}$ denote the space of
homogeneous polynomial mappings of degree $m$ of $\cala$ in
$C^n(\cala)$. The direct sum $\displaystyle{\calp = \oplusinf_{m,n}}
\calp^{m,n}$ is an associative bigraded algebra in a natural way. One
defines the total degree of an element of $\calp^{m,n}$ to be $2m+n$ ;
$\calp$ is a graded algebra for the total degree and
$C(\cala)=\displaystyle{\oplusinf_n}\calp^{0,n}$ is a graded subalgebra
of $\calp$. The composition with the differential $d$ of $C(\cala)$ is a
differential, again denoted by $d$, of the graded algebra $\calp$ which
extends the differential $d$ of $C(\cala)$ . One has
$d\calp^{m,n}\subset \calp^{m,n+1}$. By using the operation $A\mapsto
i_A$, one can define another differential, $\delta$, on $\calp$. Namely
if $\omega \in \calp$ is the polynomial mapping $A\mapsto \omega_A$ of
$\cala$ in $C(\cala)$, then $\delta\omega$ is the polynomial mapping
$A\mapsto (\delta\omega)_A=i_A\omega_A$ of $\cala$ in $C(\cala)$. One
has $\delta \calp^{m,n}\subset \calp^{m+1,n-1}$ so $\delta$ is of total
degree $2-1=1$ and the fact that $\delta$ is an antiderivation satisfying
$\delta^2=0$ follows from the fact that, for any $A\in\cala$, $i_A$ is
an antiderivation of $C(\cala)$ satisfying $i^2_A=0$. Notice that
$C^n_H(\cala)$ is the kernel of $\delta \restriction
C^n(\cala)=\calp^{0,n}$ $(:\calp^{0,n}\rightarrow \calp^{1,n-1})$.\\
{\it As a vector space}, $\calp^{m,n}$ can be identified to the
subspace of elements of $C^{m+n}(\cala)$ which are symmetric in their $m$
first arguments: For $\omega\in \calp^{m,n}$, $A\mapsto \omega_A$, there
is a unique $\xi_\omega \in C^{m+n}(\cala)$ symmetric in the $m$
first arguments such that
$$\omega_A(A_1,\dots , A_n)=\xi_\omega (\underbrace{A,\dots, A}_m,
A_1,\dots, A_n),\>  \forall A, A_i\in \cala.$$
Let $\cal I^{m,n}$ denote the subspace of $\calp^{m,n}$ consisting of the
$\omega\in \calp^{m,n}$ such that $\xi_\omega \in C^{m+n}_I(\cala)$, (i.e.
such that $\xi_\omega$ is invariant). $\cal I = \oplus {\cal I}^{m,n}$
is a graded subalgebra (also a bigraded subalgebra in the obvious sense)
of $\calp$ which is stable by $d$ and $\delta$ and, furthermore, $d$ and
$\delta$ anticommute on $\cal I$.\\
Notice that one has $\cal I^{m,0}=\cal I^m_S(\cala_{Lie})$ and $\cal
I^{0,n}=C^n_I(\cala)$ and that $C^n_B(\cala)$ is the kernel of
$\delta\restriction C^n_I(\cala)=\cal I^{0,n}$ $(:\cal I^{0,n}\rightarrow
\cal I^{1,n-1})$. The algebras $\calp$ and $\cal I$ are bigraded and $d$
and $\delta$ are bihomogeneous, therefore the $d$ and the $\delta$
cohomologies of $\calp$ and $\cal I$ are also bigraded algebras. By
using composition with the homotopy $h$ of the proof of proposition 1
and by noticing that $\cal I$ is stable by this composition, one obtains
the following generalization of proposition 1.

\begin{proposition}
One has $H^{m,n}(\calp,d)=0$, $H^{m,n}(\cal I,d)=0$ for $n\geq 1$ and
$H^{m,0}(\calp,d)=\calp^{m,0}$, $H^{m,0}(\cal I,d)=\cal
I^{m,0}=\cal I^m_S(\cala_{Lie}$).
\end{proposition}

\noindent Concerning the cohomology of $\delta$ one has the following result

\begin{proposition}
One has $H^{m,n}(\calp,\delta)=0$, $H^{m,n}(\cali,\delta)=0$ for $m\geq
1$ and $H^{0,n}(\calp,\delta)=C^n_H(\cala)$,
$H^{0,n}(\cali,\delta)=C^n_B(\cala)$.
\end{proposition}

\noindent {\bf Proof.} The last part of the proposition $(m=0)$ is
obvious since one has $H^{0,n}(\calp,\delta)=\ker(\delta\restriction
C^n(\cala))$ and $H^{0,n}(\cali,\delta)=\ker (\delta\restriction
C^n_I(\cala))$. Therefore from now on, assume that one has $m\geq 1$.
Define a linear mapping $\ell$ of $\calp$ in itself with
$\ell(\calp^{m,n})\subset \calp^{m-1,n+1}$ by $$(\ell\omega)_A
(A_1,\dots,A_{n+1})
=\frac{d}{dt}\omega_{A+tA_1}(A_2,\dots,A_{n+1})\vert_{t=0}$$ for $
\omega\in \calp^{m,n}$. One has $(\delta\ell + \ell \delta
)\omega=m\omega + \calh \omega$ where $\calh \omega$ is given by
$$(\calh\omega)_A (A_1,\dots,A_n)=
 \sum^{n+1}_{p=2} (-1)^p \omega_A (A_2,\dots,A_{p-1},
A_1, A_p,\dots,A_n),$$
$(\omega\in \calp^{m,n})$. Notice that if $\omega$
is such that $\omega_A(A_1,\dots,A_n)$ is antisymmetric in $A_1,\dots,
A_n$, then $\calh\omega=n\omega$ and therefore $\ell$ gives an homotopy
for such $\omega$. The following lemma, which is a combinatorial
statement in the algebra of the permutation group, will lead to an
homotopy for the general case. The proof of this lemma (which is
probably known) will be given in appendix.

\begin{lemma}
One has on $\calp^{m,n}, \prod^{n-2}_{p=0} (\calh - p\
id)=\prod^{n-1}_{p=0}(\calh-p\ id) =\cals$, where $\cals\omega$ is given
as before (antisymmetrisation) by $(\cals
\omega)_A(A_1,\dots,A_n)=\sum_{\pi\in\cals_n} \varepsilon(\pi)\omega_A
(A_{\pi(1)},\dots,A_{\pi(n)})$, i.e. $(\cals\omega)_A=\cals\omega_A$.

\end{lemma}

Let $\omega\in\calp^{m,n}$ with $m\geq 1$ be such that $\delta\omega=0$.
Then $\delta\ell\omega=m\omega+\calh\omega$, so one also have
$\delta\calh\omega=0$ and, by induction, $\delta\calh^p\omega=0$ for any
integer $p$; i.e. one has $\delta P(\calh)\omega=0$ for any polynomial
$P$. Define $\omega_r\in \calp^{m,n}$, for $r=1,2,\dots,n$, by
$\omega_1=\omega,\ \omega_2=\calh\omega -(n-2)\omega,\dots,\ \omega_r
=\prod^r_{p=2} (\calh-(n-p) id)\omega,\dots,\hfill\\ \omega_n=\calh(\calh -
id)\dots (\calh- (n-2) id)\omega$. One has
$\delta\ell\omega_r=m\omega_r+\calh\omega_r=\\ (m+n-r-1)\omega_r+\omega_{r+1}$,
i.e.
$$\omega_r=\delta\ell\left(\frac{\omega_r}{m+n-(r+1)}\right) -
\frac{\omega_{r+1}}{m+n-(r+1)}\vir$$ for $r\leq n-1$. This implies that
$$\omega=\delta\ell\left(\sum^{n-1}_{r=1}
\frac{(-1)^{r+1}}{\prod^{r+1}_{p=2}(m+n-p)} \omega_r\right) -
\frac{(-1)^n}{\prod^n_{p=2}(m+n-p)} \omega_n.$$ On the other hand, it
follows from the lemma and the previous discussion (antisymmetry) that
$\omega_n=\delta\ell(\frac{1}{m+n}\omega_n)$ and therefore one has an
homotopy formula, for $\omega\in\calp^{m,n}$ with $m\geq 1$ satisfying
$\delta\omega=0$, of the form $\omega=\delta\delta'\omega$ where
$\delta'=\ell\circ Q^{m,n}(\calh)$ and where the polynomial $Q^{m,n}$ is easily
computed from the previous formulae. Since $\ell$ and $\calh$ preserve
$\cali$ this achieves the proof of proposition 4.$\square$\\

The proof of the theorem 1 will now follow from $H^{m,n}(\cali,d)=0$ for
$n\geq 1$, $H^{m,n}(\cali,\delta)=0$ for $m\geq 1$,
$$H^{m,0}(\cali,d)=\cali^m_S(\cala_{Lie})\ \mbox{and}\
H^{0,n}(\cali,\delta)=C^n_B(\cala)$$
by a standard spectral sequence argument in the bicomplex
$(\cali,d,\delta)$.\\

Let $H(\delta\vert d)$ denote the
$\delta$-cohomology modulo $d$ of $\cali$, i.e. $$H^{m,n}(\delta\vert
d)=Z^{m,n}(\delta\vert d)/B^{m,n}(\delta\vert d)$$ where
$Z^{m,n}(\delta\vert d)$ is the space of the $\alpha^{m,n}\in
\cali^{m,n}$ for which there is an $\alpha^{m+1,n-2}\in \cali^{m+1,n-2}$
such that $\delta\alpha^{m,n}+d\alpha^{m+1,n-2}=0$ and where
$B^{m,n}(\delta\vert d)= \delta \cali^{m-1,n+1}+d\cali^{m,n-1}\ (\subset
\cali^{m,n})$. With these notations, one has the following result.

\begin{proposition}
One has the following isomorphisms:\\
$H^{2p}_B(\cala)\simeq
H^{k,2(p-k)-1}(\delta\vert d) \simeq \cali^p_S(\cala_{Lie})$ for $1\leq
k \leq p-2$,\\
$H^{2p+1}_B(\cala)\simeq H^{k,2(p-k)}(\delta\vert d)\simeq
0$ for $1\leq k\leq p-1$, $H^4_B(\cala)\simeq \cali^2_S(\cala_{Lie})$,
$H^3_B(\cala)\simeq 0$ and $H^2_B(\cala)\simeq \cali^1_S(\cala_{Lie})$.
\end{proposition}

\noindent {\bf Proof.}  Let $\alpha^{m,n}\in \cali^{m,n}$ be a
$\delta$-cocycle modulo $d$, i.e. there is a $\alpha^{m+1,n-2}\in
\cali^{m+1,n-2}$ such that $\delta\alpha^{m,n}+d\alpha^{m+1,n-2}=0$. By
applying $\delta$, one obtains $\delta
d\alpha^{m+1,n-2}=-d\delta\alpha^{m+1,n-2}=0$, therefore, if $n\geq 4$,
there is in view of proposition 3 a $\alpha^{m+2,n-4}\in
\cali^{m+2,n-4}$ such that $\delta\alpha^{m+1,n-2}+d\alpha^{m+2,n-4}=0$,
which means that $\alpha^{m+1,n-2}$ is also a $\delta$-cocycle modulo
$d$. If $\alpha^{m,n}$ is exact, i.e. if there are $\beta^{m-1,n+1}\in
\cali^{m-1,n+1}$ and $\beta^{m,n-1}\in \cali^{m,n-1}$ such that
$\alpha^{m,n}=\delta\beta^{m-1,n+1}+d\beta^{m,n-1}$, then
$d(\alpha^{m+1,n-2}-\delta\beta^{m,n-1})=0$ which implies, again by
proposition 3 (since $n-2\geq 2>0$), that there is a $\beta^{m+1,n-3}$
such that $\alpha^{m+1,n-2}=\delta\beta^{m,n-1}+d\beta^{m+1,n-3}$ i.e.
$\alpha^{m+1,n-2}$ is also exact. Therefore there is a well defined
linear mapping $\partial : H^{m,n}(\delta\vert d)\rightarrow
H^{m+1,n-2}(\delta\vert d)$ for $n\geq 4$ such that $\partial
[\alpha^{m,n}]=[\alpha^{m+1,n-2}]$. Let now $\alpha^{m+1,n-2}\in
\cali^{m+1,n-2}$ be a $\delta$-cocycle modulo $d$,  i.e. there is
$\alpha^{m+2,n-4}\in \cali^{m+2,n-4}$ such that
$\delta\alpha^{m+1,n-2}+d\alpha^{m+2,n-4}=0$. By applying $d$, one
obtains $\delta d\alpha^{m+1,n-2}=0$ which implies, in view of
proposition 4, that there is a $\alpha^{m,n}\in\cali^{m,n}$ such that
$\delta\alpha^{m,n}+d\alpha^{m+1,n-2}=0$. This means that $\partial$ is
surjective. Assume that $[\alpha^{m+1,n-2}]=0$ i.e.
$\alpha^{m+1,n-2}=\delta\beta^{m,n-1}+d\beta^{m+1,n-3}$ $(\beta\in\cali)$
then one has $\delta(\alpha^{m,n}-d\beta^{m,n-1})=0$ which implies that
$[\alpha^{m,n}]=0$ if $m\geq 1$ or that $\alpha^{0,n}- d\beta^{0,n-1}\in
C^n_B(\cala)$ if $m=0$, again by proposition 4. {\it Thus} $\partial
:H^{m,n}(\delta\vert d)\rightarrow H^{m+1,n-2}(\delta\vert d)$ {\it are
isomorphisms for} $n\geq 4$ {\it and} $m\geq 1$ {\it and, for} $m=0$
($n\geq 4$),
 $ \partial : H^{0,n}(\delta\vert d) \rightarrow H^{1,n-2}(\delta\vert
d)$ {\it is surjective and its kernel is the image of}
$C^n_B(\cala)=H^{0,n}(\cali,\delta)$ {\it in} $H^{0,n}(\delta\vert d)$.\\
On the other hand, if $\alpha^{0,n}\in \cali^{0,n}$ is a
$\delta$-cocycle modulo $d$, i.e. $\delta\alpha^{0,n}+d\alpha^{1,n-2}=0$
then $d\alpha^{0,n}\in C^{n+1}_I(\cala)$ is a basic cocycle of $\cala$
i.e. $d\alpha^{0,n}\in Z^{n+1}_B(\cala)$ and if $\alpha^{0,n}$ is exact,
i.e. $\alpha^{0,n}=d\beta^{0,n}$ with $\beta^{0,n}\in \cali^{0,n}$,
then $d\alpha^{0,n}=0$. Therefore, with obvious notations, one has a
linear mapping $d^\sharp:H^{0,n}(\delta\vert d)\rightarrow
H^{n+1}_B(\cala)$, $d^\sharp[\alpha^{0,n}]=[d\alpha^{0,n}]$. If $z^{n+1}\in
C^{n+1}_B(\cala)$ is closed i.e. $z^{n+1}\in Z^{n+1}_B(\cala)$ then, in
view of proposition 1, there is a $\alpha^{0,n}\in
C^n_I(\cala)=\cali^{0,n}$ such that $z^{n+1}=d\alpha^{0,n}$; one has
$d\delta\alpha^{0,n}=0$, which implies that $\alpha^{0,n}$ is a
$\delta$-cocycle modulo $d$ if $n\geq 2$ (by proposition 3). Thus $d^\sharp$
is surjective for $n\geq 2$ and one obviously has $ker(d^\sharp)$ = image of
$C^n_B(\cala)$ in $H^{0,n}(\delta\vert d)$. Applying this for  $n\geq
4$ and the previous results, one obtains isomorphisms:
$$H^{2p}_B(\cala)\simeq H^{k,2(p-k)-1}(\delta\vert d)\ \mbox{for}\
1\leq k\leq p-2$$
and
$$H^{2p+1}_B(\cala)\simeq H^{k,2(p-k)}(\delta\vert d)\ \mbox{for}\
1\leq k\leq p-1.$$
Thus, to achieve the proof, it remains to show that one has:
\begin{description}
\item[(i)] $H^{m,2}(\delta\vert d)=0$ for $m\geq 1$ and
$H^{0,2}(\delta\vert d)$ = image of $C^2_B(\cala)$
\item[(ii)] $H^{m,3}(\delta\vert d)\simeq \cali^{m+2}_S(\cala_{Lie})$ for
$m\geq 1$ and $H^{0,3}(\delta\vert d)$/image of $C^3_B(\cala)\simeq
\cali^2_S(\cala_{Lie})$
\item[(iii)] $H^2_B(\cala)\simeq \cali^1_S(\cala_{Lie})$, (remembering
that $C^1_B(\cala)=0$).
\end{description}

Let $\alpha^{m,2}\in \cali^{m,2}$ be a $\delta$-cocyle modulo $d$ then,
(since $d\alpha^{m+1,0}\equiv 0$), $\alpha^{m,2}$ is a $\delta$-cocycle,
i.e. $\delta\alpha^{m,2}=0$, which implies, by proposition 4, that
$\alpha^{m,2}\in\delta\cali^{m-1,1}$ for $m\geq 1$ and, for $m=0$,
$\alpha^{0,2}\in C^2_B(\cala)=H^{0,2}(\cali,\delta)$. This proves (i).\\
Let $\alpha^{m,3}\in\cali^{m,3}$ be a $\delta$-cocycle modulo $d$, i.e.
there is a $\alpha^{m+1,1}\in \cali^{m+1,1}$ such that
$\delta\alpha^{m,3}+d\alpha^{m+1,1}=0$. Then one has
$\delta\alpha^{m+1,1}=P^{m+2}\in
\cali^{m+2}_S(\cala_{Lie})=\cali^{m+2,0}$. If
$\alpha^{m,3}=\delta\beta^{m-1,4}+d\beta^{m,2}$ for $\beta^{m-1,4}\in
\cali^{m-1,4}$ and $\beta^{m,2}\in \cali^{m,2}$, (i.e. if $\alpha^{m,3}$
is exact), one has $d(\alpha^{m+1,1}-\delta\beta^{m,2})=0$ which
implies, by proposition 3 and by $d\cali^{m+1,0}=0$, that
$\alpha^{m+1,1}=\delta\beta^{m,2}$ and therefore
$\delta\alpha^{m+1,1}=P^{m+2}=0$. Thus there is a well defined linear
mapping $j:H^{m,3}(\delta\vert d)\rightarrow
\cali^{m+2}_S(\cala_{Lie})$, ($j([\alpha^{m,3}])=P^{m+2}$). Let
$P^{m+2}$ be an arbitrary element of $\cali^{m+2}_S(\cala_{Lie})$; then,
by proposition 4, there is a $\alpha^{m+1,1}\in\cali^{m+1,1}$ such that
$\delta\alpha^{m+1,1}=P^{m+2}$ and, since $dP^{m+2}=0$, one has $\delta
d\alpha^{m+1,1}=0$ which implies again by proposition 4 that there is a
$\alpha^{m,3}$ such that $\delta\alpha^{m,3}+d\alpha^{m+1,1}=0$. This
shows that $j$ is surjective. If $\delta\alpha^{m+1,1}=0$, then, by
proposition 4, $\alpha^{m+1,1}=\delta\beta^{m,2}$ and therefore
$\delta(\alpha^{m,3}-d\beta^{m,2})=0$ which implies again by proposition
4 that $\alpha^{m,3}=\delta\beta^{m-1,4}+d\beta^{m,2}$ if $m\geq 1$ and,
for $m=0$, $\alpha^{0,3}-d\beta^{0,2}\in C^3_B(\cala)$. This proves
(ii).\\
Finally let $z^2\in C^2_I(\cala)$ be a basic cocycle, i.e. $dz^2=0$ and
$\delta z^2=0$, then $z^2=d\alpha^1$ for a unique $\alpha^1\in
C^1_I(\cala)$ (since $dC^0(\cala)=0$ and by proposition 1). Conversely,
if $\alpha^1\in C^1_I(\cala)$ then $d\alpha^1$ is basic; therefore
$H^2_B(\cala)\simeq C^1_I(\cala)$ since $C^1_B(\cala)=0$. But one has
canonically $C^1_I(\cala)=\cali^1_S(\cala_{Lie}).\square$\\

This proves of course Theorem 1, but it is worth noticing that in the
above proof there is also a computation of the $\delta$-cohomology
modulo $d$ of $\cali$.

\section{Sketch of another approach: Connection with the Lie algebra
cohomology}

There is another way to study the basic cohomology of $\cala$ which
connects it with the Lie algebra cohomology of $\cala_{Lie}$: It is to
study the spectral sequence corresponding to the filtration of the
differential algebra $C(\cala)$ associated to the operation $i$ of the
Lie algebra $\cala_{Lie}$ in the differential algebra $C(\cala)$, [5].
This filtration $\calf$ is defined by
$$\calf^p(C^n(\cala))=\{\omega \in C^n(\cala)\vert i_{A_{1}}\dots
i_{A_{n-p+1}} (\omega)=0,\ \forall A_i\in \cala\}$$
for $0\leq p \leq n$ and $\calf^p(C(\cala))=\displaystyle{\oplusinf_{n\geq p}}
\calf^p(C^n(\cala))$.\\
One has $$\calf^0(C(\cala))=C(\cala),\  \calf^p(C(\cala)) \cdot
\calf^q(C(\cala)) \subset \calf^{p+q}(C(\cala))$$ and $$d\calf^p(C(\cala))
\subset \calf^p (C(\cala))$$ i.e. $\calf$ is a (decreasing) filtration of
graded differential algebra. To such a filtration corresponds a
convergent spectral sequence $(E_r,d_r)_{r\in \Bbb N}$, where
$E_r=\displaystyle{\oplusinf_{p,q\in \Bbb N}} E^{p,q}_r$ is a bigraded algebra
and
$d_r$ is a homogeneous differential on $E_r$ of bidegree $(r,1-r)$. The
triviality of the cohomology of $C(\cala)$, (i.e. proposition 1),
implies that $E^{p,q}_\infty=0$ for $(p,q)\not=(0,0)$ and
$E^{0,0}_\infty= \k$. The spectral sequence starts with the graded space
$E_0$ associated to the filtration i.e.
$E^{p,q}_0=\calf^p(C^{p+q}(\cala))/\calf^{p+1}(C^{p+q}(\cala))$ and $d_0$
is induced by the differential $d$ of $C(\cala)$. If
$\omega\in\calf^p(C^{p+q}(\cala))$ then $i_{A_{1}}\dots i_{A_{q}}\omega$
is in $C^p_H(\cala)$ and is antisymmetric in $A_1,\dots,A_q$. Therefore
$(A_1,\dots,A_q)\mapsto i_{A_{1}}\dots i_{A_{q}}\omega$ is a $q$-cochain of
the Lie algebra $\cala_{Lie}$ with values in $C^p_H(\cala)$ for the
representation $A\mapsto L_A$ of the Lie algebra $\cala_{Lie}$ in
$C^p_H(\cala)$. This defines a linear map of $\calf^p(C^{p+q}(\cala))$
in the space of $q$-cochains of $\cala_{Lie}$ with values in
$C^p_H(\cala)$. The kernel of this map is, by definition,
$\calf^{p+1}(C^{p+q}(\cala))$. In our case, it is straightforward to
show that this map is surjective, i.e. that $E^{p,q}_0$ identifies with the
space of $q$-cochains of the Lie algebra $\cala_{Lie}$ with values in the
space $C^p_H(\cala)$ of horizontal elements of $C^p(\cala)$ and that
then, $d_0$ coincides with the Chevalley-Eilenberg differential. Thus
$E_1=H(E_0,d_0)$ is the Lie algebra cohomology of $\cala_{Lie}$ with
value in $C_H(\cala)$, $E^{p,q}_1=H^q(\cala_{Lie}, C^p_H(\cala))$. In
particular $E^{0,\ast}_1$ is the ordinary cohomology of $\cala_{Lie}$
(i.e. with value in the trivial representation in $\k$) and
$E^{\ast,0}_1$ is the space of invariant elements of $C_H(\cala)$, i.e.
the space $C_B(\cala)$ of basic elements of $C(\cala)$,
$E^{n,0}_1=C^n_B(\cala)$. Furthermore, on $E^{\ast,0}_1=C_B(\cala)$,
$d_1$ is just the differential $d$ of $C(\cala)$ restricted to
$C_B(\cala)$. Therefore $E^{\ast,0}_2$ is the basic cohomology
$H_B(\cala)$ of $\cala$, $E^{n,0}_2=H^n_B(\cala)$. This shows that the
spectral sequence connects the basic cohomology of $\cala$ to the Lie
algebra cohomology of its underlying Lie algebra $\cala_{Lie}$. The
connection between the Lie algebra cohomology of $\cala_{Lie}$ and the
ad$^\ast$-invariant polynomials, i.e. $H_B(\cala)$ in our case, is well
known but an interest of the last approach could be to catch the
primitive parts.

\section*{Appendix: Proof of Lemma 1}

Let $\cals_n$ be the group of permutations of $\{1,\dots,n\}$.\\
In the algebra of this group, let us define the antisymmetrisation
operator
$$\cals =\sum_{\pi\in\cals_n}\varepsilon(\pi)\pi$$
and the operators
$$\calh_{(k)}=\sum_{{\pi\in\cals_n}\atop
{\pi^{-1}(k+1)<\dots<\pi^{-1}(n)}}\varepsilon(\pi)\pi$$
for any $1\leq k\leq n$, where $\varepsilon(\pi)$ denotes the signature
of the permutation $\pi$.\\
Notice that $$\calh_{(n)}=\calh_{(n-1)}=\cals$$
and one easily shows that
$$\calh\equiv \calh_{(1)} = \sum^n_{p=1}(-1)^{p+1}\gamma_p$$
where $\gamma_p$ is the permutation $(1,\dots,p,\dots,n)\mapsto
(2,\dots,p,1,p+1,\dots,n)$.\\
With these definitions, one has the following result:\\

\noindent{\bf LEMMA}
{\it For any} $1\leq k\leq n-1$,
$$\calh\calh_{(k)} = k \calh_{(k)}+ \calh_{(k+1)}$$
{\bf Proof.}\\
\begin{eqnarray*}
\calh\calh_{(k)} & = &\left(\sum^n_{p=1}(-1)^{p+1}\gamma_p\right) \left(
\sum_{{\pi\in\cals_n}\atop{\pi^{-1}(k+1)<\dots<\pi^{-1}(n)}}
\varepsilon(\pi)\pi\right)
\\
& = & \sum^k_{p=1}
\sum_{{\pi\in\cals_n}\atop{\pi^{-1}(k+1)<\dots<\pi^{-1}(n)}}(-1)^{p+1}\varepsilon(\pi)\gamma_p\pi\\
& & \mbox{}+  \sum^n_{p=k+1} \sum_{{\pi\in\cals_n}\atop{\pi^{-1}(k+1)<\dots
<\pi^{-1}(n)}}
(-1)^{p+1}\varepsilon(\pi)\gamma_p\pi
\end{eqnarray*}
Now, define $\pi'=\gamma_p\pi\in \cals_n$; one has
$\varepsilon(\pi')=(-1)^{p+1}\varepsilon(\pi)$.\\
For $p\leq k$, one has $\pi'^{-1}(q)=\pi^{-1}(q)$ for any $k+1\leq q\leq
n$.\\
So, in the first summation, for a fixed $p$, the sum
over the $\pi\in\cals_n$ such that $\pi^{-1}(k+1)<\dots<\pi^{-1}(n)$ can be
replaced
by the sum over the $\pi'\in\cals_n$ such that
$\pi'^{-1}(k+1)<\dots<\pi'^{-1}(n)$.
Thus
$$\sum^k_{p=1} \sum_{{\pi\in\cals_n}\atop{\pi^{-1}(k+1)<\dots<\pi^{-1}(n)}}
(-1)^{p+1}\varepsilon(\pi)\gamma_p\pi=\sum^k_{p=1}
\sum_{{\pi'\in\cals_n}\atop{\pi'^{-1}(k+1)<\dots<\pi'^{-1}(n)}}
\varepsilon(\pi')\pi'=k\calh_{(k)}$$
Now, for $p\geq k+1$, one has $\pi'^{-1}(q)=\pi^{-1}(q-1)$ for any
$k+2\leq q \leq p$ and $\pi'^{-1}(q)=\pi^{-1}(q)$ for any $p+1\leq q
\leq n$. So one has only
$$\pi'^{-1}(k+2)<\dots<\pi'^{-1}(n),$$
and in the second summation the sum over $p$ and $\pi$ can be replaced by the
sum over the
$\pi'\in \cals_n$ such that $\pi'^{-1}(k+2)<\dots<\pi'^{-1}(n).$
Thus
$$\sum^n_{p=k+1}\sum_{{\pi\in\cals_n}\atop{\pi^{-1}(k+1)<\dots<\pi^{-1}(n)}}
(-1)^{p+1}\varepsilon(\pi)\gamma_p\pi = \calh_{(k+1)}.\  \square$$

By induction, this lemma shows that for any $1\leq k\leq n$
$$\calh_{(k)} = \prod^{k-1}_{p=0} (\calh - p\ id)$$
where we recall $\calh\equiv \calh_{(1)}$.\\

So for $k=n$ and $k=n-1$, one has
$$\calh_{(n)} =  \prod^{n-1}_{p=0}(\calh- p\ id) = \cals$$
$$\calh_{(n-1)} = \prod^{n-2}_{p=0}(\calh- p\ id) = \cals$$
Now, notice that the operators $\calh$ and $\cals$ of Lemma 1 are
representations of the operators $\calh$ and $\cals$ above in the linear
space $\calp^{m,n}$ (in fact only in $C^n(\cala))$. This proves Lemma 1.

\end{document}